%% file: 2014arxiv.tex
\documentclass[runningheads]{llncs}
\usepackage{llncsdoc}
\usepackage{amssymb,amsmath}
\usepackage{graphicx}
\usepackage{Lilypack}
\usepackage{url}
\usepackage{subfigure}
\usepackage{multirow}
\usepackage{tablefootnote}
\urldef{\mailsa}\path|{ioanniskon,ratliffl,jinming,sastry,spanos}@eecs.berkeley.edu|
\begin{document}

\mainmatter  % start of an individual contribution

\title{Social Game for Building Energy Efficiency:\\ Utility
Learning, Simulation, and Analysis}

\titlerunning{Social Game for Energy Efficiency}

\author{Ioannis C. Konstantakopoulos$^{\star\star}$ \and Lillian J.
Ratliff\thanks{This work is supported by NSF
CPS:Large:ActionWebs award number 0931843, TRUST (Team for Research in
Ubiquitous Secure Technology) which receives support from NSF (award
number CCF-0424422), and FORCES (Foundations Of Resilient
CybEr-physical Systems) which receives support from NSF (award number
CNS-1239166)} \and Ming Jin$^{\star\star}$\and S. Shankar
Sastry$^\star$\and Costas Spanos\thanks{This research is funded by the Republic of Singapore's National Research Foundation through a grant to the Berkeley Education Alliance for Research in Singapore (BEARS) for the Singapore-Berkeley Building Efficiency and Sustainability in the Tropics (SinBerBEST) Program. BEARS has been established by the University of California, Berkeley as a center for intellectual excellence in research and education in Singapore.}}%
%\thanks{Please note that the LNCS Editorial assumes that all authors have used
%the western naming convention, with given names preceding surnames. This determines
%the structure of the names in the running heads and the author index.}%
%\and Ursula Barth\and Ingrid Haas\and Frank Holzwarth\and\\
%Anna Kramer\and Leonie Kunz\and Christine Rei\ss\and\\
%Nicole Sator\and Erika Siebert-Cole\and Peter Stra\ss er}

\institute{Department of Electrical Engineering and Computer Sciences\\
University of California, Berkeley,\\
Berkeley, California, 94720\\
\mailsa}
%\url{http://www.springer.com/lncs}}


%\toctitle{Lecture Notes in Computer Science}
%\tocauthor{Authors' Instructions}
\maketitle

\begin{abstract}
We describe a social game that we designed for encouraging energy efficient
behavior amongst building occupants with the aim of reducing overall energy
consumption in the building.  Occupants vote for their desired lighting level
and win points which are used in a lottery based on how far their vote is from
the maximum setting. We assume that the
occupants are utility maximizers and that their utility functions capture the
tradeoff between winning points and their comfort level. We model the occupants
as
non-cooperative agents in a continuous game and we characterize their play using
the Nash
equilibrium concept. Using occupant voting data, we parameterize their utility functions and use a convex
optimization problem to estimate the parameters. We simulate the game defined by
the estimated utility functions and show that the estimated model for occupant
behavior is a good predictor of their actual behavior. In addition, we show that
due to the social game, there is a significant reduction in energy
consumption.
%resulting
%from the estimated parameters 
  \keywords{Utility Learning, Energy Efficiency, Game Theory}
\end{abstract}

%%%%%%%%%%%%%%%%%%%%%%%%%%%%%%%%%%%%%%%%%%%%%%%%%%%%%%%%%%%%%%%%%%%%%%%%%%%%%%%
%%%%%%%%%%%%%%%%%%%%%%%%%%%%%%%%%%%%%%%%%%%%%%%%%%%%%%%%%%%%%%%%%%%%%%%%%%%%%%%
%%%%%%%%%%%%%%%%%%%%%%%%%%%%%%%%%%%%%%%%%%%%%%%%%%%%%%%%%%%%%%%%%%%%%%%%%%%%%%%
\section{Introduction}
\input{intro}

%%%%%%%%%%%%%%%%%%%%%%%%%%%%%%%%%%%%%%%%%%%%%%%%%%%%%%%%%%%%%%%%%%%%%%%%%%%%%%%
%%%%%%%%%%%%%%%%%%%%%%%%%%%%%%%%%%%%%%%%%%%%%%%%%%%%%%%%%%%%%%%%%%%%%%%%%%%%%%%
%%%%%%%%%%%%%%%%%%%%%%%%%%%%%%%%%%%%%%%%%%%%%%%%%%%%%%%%%%%%%%%%%%%%%%%%%%%%%%%
%\section{Background}
%\label{sec:back}
%\input{back}
%%%%%%%%%%%%%%%%%%%%%%%%%%%%%%%%%%%%%%%%%%%%%%%%%%%%%%%%%%%%%%%%%%%%%%%%%%%%%%%
%%%%%%%%%%%%%%%%%%%%%%%%%%%%%%%%%%%%%%%%%%%%%%%%%%%%%%%%%%%%%%%%%%%%%%%%%%%%%%%
%%%%%%%%%%%%%%%%%%%%%%%%%%%%%%%%%%%%%%%%%%%%%%%%%%%%%%%%%%%%%%%%%%%%%%%%%%%%%%%

%%%%%%%%%%%%%%%%%%%%%%%%%%%%%%%%%%%%%%%%%%%%%%%%%%%%%%%%%%%%%%%%%%%%%%%%%%%%%%%
%%%%%%%%%%%%%%%%%%%%%%%%%%%%%%%%%%%%%%%%%%%%%%%%%%%%%%%%%%%%%%%%%%%%%%%%%%%%%%%
%%%%%%%%%%%%%%%%%%%%%%%%%%%%%%%%%%%%%%%%%%%%%%%%%%%%%%%%%%%%%%%%%%%%%%%%%%%%%%%
\section{Game Formulation}
\label{sec:game}
\input{game}

\section{Experimental Set-Up}
\label{sec:experiment}
\input{experiment}

\section{Results}
\label{sec:results}
\input{results}

%%%%%%%%%%%%%%%%%%%%%%%%%%%%%%%%%%%%%%%%%%%%%%%%%%%%%%%%%%%%%%%%%%%%%%%%%%%%%%%
\section{Discussion and Future Work}
\label{sec:discussion}
\input{discussion}

\bibliographystyle{splncs03}
\bibliography{2014arxiv}
%\begin{thebibliography}{4}
%
%\bibitem{jour} Smith, T.F., Waterman, M.S.: Identification of Common Molecular
%Subsequences. J. Mol. Biol. 147, 195--197 (1981)
%
%\bibitem{lncschap} May, P., Ehrlich, H.C., Steinke, T.: ZIB Structure Prediction Pipeline:
%Composing a Complex Biological Workflow through Web Services. In: Nagel,
%W.E., Walter, W.V., Lehner, W. (eds.) Euro-Par 2006. LNCS, vol. 4128,
%pp. 1148--1158. Springer, Heidelberg (2006)
%
%\bibitem{book} Foster, I., Kesselman, C.: The Grid: Blueprint for a New Computing
%Infrastructure. Morgan Kaufmann, San Francisco (1999)
%
%\bibitem{proceeding1} Czajkowski, K., Fitzgerald, S., Foster, I., Kesselman, C.: Grid
%Information Services for Distributed Resource Sharing. In: 10th IEEE
%International Symposium on High Performance Distributed Computing, pp.
%181--184. IEEE Press, New York (2001)
%
%\bibitem{proceeding2} Foster, I., Kesselman, C., Nick, J., Tuecke, S.: The Physiology of the
%Grid: an Open Grid Services Architecture for Distributed Systems
%Integration. Technical report, Global Grid Forum (2002)
%
%\bibitem{url} National Center for Biotechnology Information, \url{http://www.ncbi.nlm.nih.gov}
%
%\end{thebibliography}

%\section*{Appendix: Springer-Author Discount}

\end{document}

%% file: intro.tex
% intro file

 Energy consumption of buildings, both residential and commercial, accounts for
 approximately 40\% of all energy usage in the
 U.S.~\cite{mcquade2009}. Lighting is a major consumer of energy in commercial
 buildings; one-fifth of all energy consumed in buildings is due to
 lighting~\cite{eia:2009ab}.
 
There have been many approaches to improve energy efficiency of buildings
through control and automation as well as incentives and pricing. From the meter
 to the consumer, many control methods, such as model predictive control, have
 been proposed as a means to improve the efficiency of building operations (see, e.g.,
 \cite{Aswani:2012kx},\cite{boman:1998aa},\cite{bourgeois:2006aa},\cite{maAnderson2011ACC},\cite{oldewurtel2010ACC},\cite{lovett:2013aa}). From
 the meter to the energy utility, many economic solutions have been proposed,
 such as dynamic pricing and smart meter technology, to reduce consumption by
 providing economic incentives (see, e.g.,\cite{Mathieu2012},\cite{Dahleh2010smartCom}).

 Many of the past approaches to building energy management only focus on heating
 and cooling of the building. We are advocating that due to new technological
 advances in building automation, incentives can be designed around more than
 just heating, ventilation and air conditioning (HVAC) systems. In particular,
 our experimental set-up allows us to design incentives based on 
 lighting and individual plug-load in addition to HVAC and interact with
 occupants through a social game.

 Social games have been used to encourage energy efficient behavior in
 transportation~\cite{merugu:2009aa} as well as in the healthcare domain for
 understanding the tradeoff between privacy and desire to win by expending
 calories~\cite{bestick:2013aa}. 

 There are many ways in which a building manager can be motivated to
 encourage energy efficient behavior. The most obvious is that they pay the
 bill or, due to some operational excellence measure, are required to maintain
 an energy effcient building. Beyond these motivations, recently demand response
 programs have begun to be implemented by utility
 companies with the goal of correcting for improper load forecasting (see, e.g.,
 \cite{albadi:2008aa},\cite{mathieu:2011aa}, \cite{m.-lee:2013aa}). In such a
 program, consumers enter into a contract with the utiltiy company in which they
 agree to change their demand in accordance with some agreed upon schedule. In
 this scenario, the building manager may now be required to keep this schedule.

Our approach to efficient building energy management focuses on office
 buildings and utilizes new builidng automation products such as the Lutron
 lighting system\footnote{\tt http://www.lutron.com/en-US/Pages/default.aspx}.
 We design a social game aimed at incentivizing occupants to
 modify their behavior so that the overall energy consumption in the building is
 reduced. 
The social game consists of occupants logging their vote for the
 lighting setting in the office and they win points based on how energy
 efficient their vote is compared to other occupants. The average of the votes
 is what is actually implemented in the office. The points are used
to determine an occupants likelihood of winning in a lottery.
 %Occupants can exchange points for
 %a chance to win in a lottery.
We designed an online platform so that occupants can login and vote, view their
points, and observe all occupants consumption patterns and points. This platform
also store all the past data allowing us to use it for estimation of the
behavior of the occupants.
%for a social game in which building occupants can vote on the
%desired lighting level above their cubicle. 

%The occupants win points
%based on how far their vote is from the baseline setting, i.e. the typical setting that
%office operated at prior to the implementation of the game. 

In this paper we present the results of a social game focused only on the
encouraging more energy efficient lighting usage; however, we emphasize that the
framework is easily adapted to incorporate the full capabilities of the
automation installed in our experimental set-up (i.e. lighting, HVAC, and plug-load).
The occupants are modeled as utility maximizers who engage in a non-cooperative game
with all occupants.  We parameterize their utility
functions in such a way that we capture the tradeoff between the desire to win
and comfort. Using data from the social game that occured over the period
of roughly three months, we formulate the utility learning problem as a convex
optimization problem and form estimates of each occupants utility function. We simulate the game using the estimated utility functions
and show that the Nash equilibrium from the simulations is a good predictor of
occupant behavior. Our results are compared to other estimation techniques. 

A major advantage of modeling occupants as utility maximizers competing
in a game and using the Nash equilibrium concept is this game theoretic model
fits in the Stackelberg framework for incentive design in which the builiding
manager performs an online estimation of occupant's utility function and designs
incentives for behavior modification. This, in essence, is a problem of
closing-the-loop around the occupants so that the building manager achieves
sustained energy savings. We leave this as future work.
%We remark that our method of modeling occupants as utility maximizers competing
%in a game and using the Nash equilibrium concept is powerful because we can
%consider the incentive design problem as a Stackelberg game in which the
%building manager estimates the utility functions and designs the incentives in
%an online framework. 
%\textcolor{red}{Lily: This needs some work, we also need to put in some more
%background on other approaches to incentive design}.
%The social game for energy savings that we have designed is such that 
%occupants in an office builing vote according to their usage preferences of shared resources and are
%rewarded with points based on how \emph{energy efficient} their strategy is in
%comparison with the other occupants. Having points increases the likelihood of
%the occupant winning in a weekly lottery. 

The rest of the paper is organized as follows.
%In Section~\ref{sec:back} we
%discuss relevant literature on similar social game experiments. 
In Section~\ref{sec:game}, we start with the game theoretic framework for
modeling the competitive environment between the non-cooperative occupants. We
formulate the utility estimation problem as a convex optimization problem and
take a dynamical systems perspective for developing a method of computing the
Nash equilibrium of the estimated game. In
Sections~\ref{sec:experiment} and \ref{sec:results} we describe the experimental
set-up and report on our findings including utility estimation results as well
as simulation of the game corresponding to the estimated utilities. Finally, in
Section~\ref{sec:discussion}, we make concluding remarks and comment on future
research directions.

%% file: game.tex
We begin by describing the game theoretic framework used for modeling the
interaction between the occupants. We remark that the use of game theory for
modeling the behavior of the occupants has several advantages. First, it is a
natural way to model agents competiting over scarce resources. It can also be
leveraged in the design of incentives for behavioral change in that it
incorporates the ability to model the occupants as strategic players.

Let the number of occupants participating in the game be denoted by $n$. 
We model the occupants as utility maximizers having utility functions composed
of two terms that capture the tradeoff between comfort and desire to win. We
model their comfort level using a Taguchi loss function which is interpreted as
modeling occupant dissatisfaction as increasing as variation increases from
their desired lighting setting. In particular, each occupant has the following
Taguchi loss function as one component of their utility function:
\begin{equation}
  \psi_i(x_i,x_{-i})=-\left( \bar{x}-x_i \right)^2
  \label{eq:taguchi}
\end{equation}
where $x_i\in\mb{R}$ is occupant $i$'s lighting vote, $x_{-i}=\{x_1, \ldots, x_{i-1},
x_{i+1}, \ldots, x_n\}$, and
\begin{equation}
  \bar{x}=\frac{1}{n}\sum_{i=1}^nx_i
  \label{eq:avg}
\end{equation}
is the average of all the occupant votes and is the lighting setting which is
implemented.

Each occupant's desire to win is modeled using the following function
\begin{equation}
  \phi_i(x_i,x_{-i})=\ln\left( \rho\frac{x_b-x_i}{nx_b-\sum_{j=1}^nx_j} \right)
  \label{eq:points}
\end{equation}
where $\rho$ is the total number of points distributed by the building manager and $x_b$ is
the baseline setting for the lights. The term inside the natural log function is
how the points are distributed; $\rho$, being the total number of points, is
multiplied by the distance an occupant's vote is from the baseline and then
normalized by the sum of the differences of all occupants' votes from the
baseline.

Hence, each occupant's utility function is given by
\begin{equation}
  f_i(x_i, x_{-i})=\psi_i(x_i, x_{-i})+\theta_i\phi_i(x_i, x_{-i})
  \label{eq:utility}
\end{equation}
where $\theta_i$ is an unknown parameter.

The occupants face the following optimization problem
\begin{equation}
\max\limits_{x_i\in S_i} f_i(x_i, x_{-i})
 %& \mathrm{s.t.}\  x_i\in [0,100]
%\end{array}\right.
  \label{eq:P1}
\end{equation}
where $S_i=[0,100]\subset\mb{R}$ is the constraint set for each $x_i$. 

Note that
each occupant's optimization problem is dependent on the other occupant's choice
variables.

We can explicitly write out the constraint set as follows. Let
$h_{i,j}(x_i, x_{-i})$ for $j\in\{1, 2\}$ denote the constraints
 on occupant $i$'s optimization problem. In particular, following
 Rosen~\cite{Rosen:1965tg}, for occupant $i$, the constraints are
 \begin{align}
   h_{i,1}(x_i)&=100-x_i\\
   h_{i,2}(x_i)&=x_i
   \label{eq:constraints}
 \end{align}
 so that we can define $\mc{C}_i=\{x_i\in \mb{R}|\ h_{i,j}(x_i)\geq 0, \ j\in\{1,2\}\}$ and
 $\mc{C}=\mc{C}_1\times\cdots\times \mc{C}_n$.
%We define the joint
%strategy space to be 
Thus, the occupants are non-cooperative agents in a continuous game with
constraints. 
%We denote the joint strategy space by $\mc{C}=S_1\times \cdots
%\times S_n\subset \mb{R}^n$. 
We model their interaction using the Nash equilibrium concept.
\begin{definition}
  A point $x\in \mc{C}$ is a {\bf Nash equilibrium} for the game $(f_1, \ldots,
  f_n)$ on $\mc{C}$
if
\begin{equation}
  f_i(x_i, x_{-i})\geq f_i(x_i',x_{-i})\ \ \forall \ x_i'\in S_i
  \label{eq:ineq}
\end{equation}
for each $i\in \{1, \ldots, n\}$. 
\end{definition}
The interpretation of the definition of Nash is as follows: no player can
unilaterally deviate and increase their cost.

If the parameters $\theta_i\geq0$, then the game is a concave
$n$-person game on a convex set. 
\begin{theorem}[\cite{Rosen:1965tg}]
  A Nash equilibrium exists for every concave $n$-person game.
  \label{thm:existence}
\end{theorem}
%Theorem 1 in the seminal work by Rosen guarantees that a Nash
%equilibrium exists. 
Define the Lagrangian of each players optimization problem as follows:
\begin{equation}
  L_i(x_i, x_{-i}, \mu_i)=f_i(x_i, x_{-i})+\sum_{j\in
  A_i(x_i)}\mu_{i,j}h_{i,j}(x_i)
  \label{eq:lagrangian}
\end{equation}
where $A_i(x_i)$ is the active constraint set at $x_i$. 
We can define 
\begin{equation}
  \omega(x)=\bmat{D_1L_1(x, \mu_i)\\ \vdots\\ D_nL_n(x, \mu_i)}
  %+\text{diag}(D_1h_1, \ldots,
  %D_nh_n)\bmat{\lambda_1& \cdots &\lambda_n}
  \label{eq:omega}
\end{equation}
where $D_iL_i$ denoets the derivative of $L_i$ with respect to $x_i$.
%performs the action of taking the derivative of
%its operand with respect to the $i$-th component of its argument. 

It is the local representation of the differential game
form~\cite{ratliff:2013aa} corresponding to the game between the occupants.

%\begin{definition}
%  A function $g(x)$ is said to be {\bf diagonally strictly concave} for $x\in
%  \mb{C}$ if for every $y,z\in \mc{C}$ 
%  \begin{equation}
%    (y-z)^Tg(z)+(z-y)^Tg(y)>0
%    \label{eq:dsc}
%  \end{equation}
%\end{definition}
\begin{definition}[\cite{ratliff:2013aa}]
  A point $x^\ast\in{\mc{C}}$ is a {\bf differential Nash equilibrium} for the game
  $(f_1, \ldots, f_n)$ on $\mc{C}$ if $\omega(x^\ast, \mu_i^\ast)=0$ and
  $D_{ii}L_i(x^\ast, \mu_i^\ast)<0$ where $\mu_{i,j}\geq 0$ for $j\in
  A_i(x_i^\ast)$.
  %or if $x^\ast\in\partial \mc{C}$, 
  %\begin{equation}
  %  D_if_i(x^\ast)+\sum_{j\in A_i(x^\ast)} \lambda_i^j D_ih_i^j(x^\ast)=0\ \ \text{for each} \
  %  i\in\{1, \ldots, n\}
  %  \label{eq:lag}
  %\end{equation}
  %where $A_i(x^\ast)$ is the set of active constraints for player $i$ at $x^\ast$ and
  %$\lambda_i^j$ is the Lagrange multiplier corresponding to constriant $h_i^j$ of player $i$'s optimization problem
  %over $x_i$ with $x_{-i}^\ast$ fixed.
\end{definition}
A sufficient condition guaranteeing that a Nash equilibrium $x$ is isolated is
that the Jacobian of $\omega(x)$, denoted
$D\omega(x)$, is invertible~\cite{ratliff:2013aa},\cite{Rosen:1965tg}.

\subsection{Utility Estimation}
We formulate the utility estimation problem as a convex optimization
problem by using first-order
necessary conditions for Nash equilibria. In particular,
the gradient of each occupant's utility function should be identically zero at
the observed Nash equilibrium. This is the case since
the observed Nash equilibria are all inside the feasible region so that none of
the constraints are active, i.e. we do not have to check the derivative of
Lagrangian of each occupant's optimization problem.

In particular, for each observation $x^{(k)}$, we assume that it corresponds to
occupants playing a strategy that is approximately a Nash
equilibrium where the superscript notation $(\cdot)^{(k)}$ indicates the $k$-th observation.
\begin{definition}
  A point $x\in \mc{C}$ is a {\bf $\varepsilon$-Nash equilibrium} for the game $(f_1, \ldots,
  f_n)$ on $\mc{C}$
if
\begin{equation}
  f_i(x_i, x_{-i})\geq f_i(x_i',x_{-i})-\vep\ \ \forall \ x_i'\in S_i
  \label{eq:ineq}
\end{equation}
for each $i\in \{1, \ldots, n\}$.  
\end{definition}

Thus, we can consider first-order optimality conditions for each
occupants optimization problem and define a residual function capturing the
amount of
sub-optimality of each occupants choice $x_i^{(k)}$~\cite{keshavarz:2011aa},\cite{ratliff:2014aa}.
Note that all our observations are on the interior of
the constraint set so we need only consider the following residual defined by
the stationarity condition for each occupant's optimization problem:
\begin{equation}
  r_{i}^{(k)}(\theta_i)=D_if_i(x_i^{(k)}, x_{-i}^{(k)})=D_i\psi_i(x_i^{(k)},
  x_{-i}^{(k)})+\theta_iD_i\phi_i(x_i^{(k)}, x_{-i}^{(k)})
  \label{eq:dzero}
\end{equation}
Define $r^{(k)}(\theta)=[r_1^{(k)}(\theta_1)\ \cdots \ r_n^{(k)}(\theta_n)]^T$.

Given observations $\{x^{(k)}\}_{k=1}^K$ where each $x^{(k)}\in \mc{C}$, we can
solve the following convex optimization problem:
\begin{equation}
  \min\limits_{\theta}\left\{\sum_{k=1}^K\chi(r^{(k)}(\theta))\bigg|\ \theta_i\geq
  0\ \ \forall \ i\{1, \ldots, n\}\right\}
  \label{eq:penaltyprob}
\end{equation}
where $\chi:\mb{R}^n\rar \mb{R}_+$ is a nonnegative, convex penalty function
satisfying $\chi(z)=0$ if and only if $z=0$, i.e. any norm on $\mb{R}^n$.
%Suppose we have $K$ observations. 

With a specific choice of $\chi$ we can explicitly write the estimation problem
as follows. Let 
\begin{equation}
  \Psi_i=\bmat{D_i\psi_i(x_i^{(1)},
  x_{-i}^{(1)})\\ \vdots\\ D_i\psi_i(x_i^{(K)},
  x_{-i}^{(K)})}, \ \Phi_i=\bmat{D_i\phi_i(x_i^{(1)},
  x_{-i}^{(1)})\\ \vdots\\ D_i\phi_i(x_i^{(K)},
  x_{-i}^{(K)})} 
  \label{eq:caps}
\end{equation}
for each $i\in\{1, \ldots, n\}$ and denote $\theta=[\theta_1\ \cdots \theta_n]^T$.
Then, we can formulate the following convex optimization problem to solve for $\theta$:
\begin{align}
  \min\limits_{\theta} \left\{\sum_{i=1}^n\|\Psi_i+\theta_i\Phi_i\|_2^2\ \bigg|\
  \theta_i\geq 0 \ \  \forall \ i\in\{1, \ldots, n\}\right\}
  %&\text{s.t.}\ \theta_i\geq 0, \ \forall \ i\in\{1, \ldots, n\}\notag
  \label{eq:est}
\end{align}
Note the constraint that the $\theta_i$'s be non-negative. This is to ensure
that the estimated utility functions are concave. We add this restriction so
that we can employ techniques from simulation of dynamical systems to the
computation of the Nash equilibrium in the resulting $n$-person concave game
with convex constraints.

\subsection{Dynamical Systems Perspective}
We can take a dynamical systems perspective in order to come up with a method
for computation of the Nash equilibrium (see, e.g. \cite{flam:1990aa},
\cite{ratliff:2013aa}, \cite{Rosen:1965tg}). We first write down a reasonable
set of dynamics, then we show that a Nash equilibrium is a stable fixed point of
these dynamics, and finally we suggest a subgradient projection
method for computation.

It is natural to consider computing Nash equilibria by following the
gradient of each occupant's utility function. Hence, we consider the dynamical
system obtained by taking the derivative with respect to their choice variable
of the Largrangian's for each occupant's optimization problem. 

%Indeed, let
%$h_{i,j}(x_i, x_{-i})$ for $j\in\{1, 2\}$ denote the constraints
% on occupant $i$'s optimization problem. In particular, following
% Rosen~\cite{Rosen:1965tg}, for occupant $i$, the constraints are
% \begin {align}
%   h_{i,1}(x_i)&=100-x_i\\
%   h_{i,2}(x_i)&=x_i
%   \label{eq:constraints}
% \end{align}
% so that we can define $\mc{C}_i=\{x_i\in \mb{R}|\ h_{i,j}(x_i)\geq 0, \ j\in\{1,2\}\}$ and
% $\mc{C}=\mc{C}_1\times\cdots\times \mc{C}_n$. 
 Due to the fact that our constraint set
 is convex, closed and bounded in $\mb{R}^n$ and there is a point in its strict interior, we
 satisfy a constraint qualification condition which is a sufficient condition
 for the Karush-Khun-Tucker (KKT) conditions for each occupant's optimization
 problem~\cite{arrow:1961aa}. It is known that for concave games, i.e. concave
 player utility functions constrained on a convex set, given that the problem
 satisfies a constraint qualification condition, then a point satisfying KKT
 conditions for each player's optimization problem is a Nash
 equilibrium~\cite{Rosen:1965tg}.
 
% In addition, it is clear that each $h_{i,j}$ is smooth. 
%Thus, the KKT conditions are given as follows. Let $x^\ast=(x_1^\ast, \ldots,
%x_n^\ast)$ be a Nash equilibrium. Then, $h_{i,j}(x_i^\ast)\geq 0$ for each
%$i\in\{1, \ldots, n\}$ and $j\in \{1, 2\}$. Further, for each
%$i\in\{1, \ldots, n\}$ there exists $\mu_{i,j}^\ast\geq 0$, for $j\in\{1, 2\}$ such that 
%$\mu_{i,j}^\ast h_{i,j}(x_i^\ast)=0$ and 
%\begin{equation}
%  0=D_if_i(x^\ast)+\mu_{i,1}^\ast D_ih_{i,1}(x_i^\ast)+\mu_{i,2}^\ast
%  D_ih_{i,2}(x_i^\ast).
%  \label{eq:kkt1}
%\end{equation}
%We remark that the KKT conditions are necessary for optimality of each
%occupant's individual optimization problem and for a Nash equilibrium. 

 We can study the continuous-time dynamical system generated by the gradient of
 the Lagrangian of each occupant's optimization problem with respect to her own
 choice variable; we let
\begin{equation}
  \dot x_i = D_if_i(x_i, x_{-i})+\sum_{j=1}^{2}\mu_{i,j}D_ih_{i,j}(x_i)
  \label{eq:dynamics}
\end{equation}
for $i\in \{1, \ldots, n\}$ and where $\mu_{i,j}$ is the $j$-th dual variable for
occupant $i$'s optimization problem. The first term is the derivative of
occupant $i$'s utility with respect to her own choice variable $x_i$. The second
term, with the appropriate dual variables $\mu_{i,j}$, ensures that for any
initial condition in the feasible set $\mc{C}$, the trajectory solving
\eqref{eq:dynamics} remains in $\mc{C}$. The right-hand side of
\eqref{eq:dynamics} is the projection of the psuedogradient on the manifold
formed by the active constraints at $x$~\cite{Rosen:1965tg}.

We can rewrite the dynamics in a compact form as follows. Let $H(x)=[Dh_1\
Dh_2]$ where $h_j(x)=[h_{1,j}\ \cdots \ h_{n,j}]^T$ for $j\in\{1,2\}$
%\begin{equation}
% h_j=\bmat{h_{1,j}\\ \vdots\\ h_{n,j}}\ \ \forall j\in\{1,2\}
%  \label{eq:h_j}
%\end{equation}
and $D$ is the Jacobian operator. Also, let $\mu=[\mu_{1, 1}\ \cdots \
\mu_{n,1}\ \mu_{1,2}\ \cdots\ \mu_{n,2}]^T$. Define $F(x,\mu)=\omega(x)+H(x)\mu
$.
%\begin{equation}
%  F(x,\mu)=\omega(x)+H(x)\mu
%  \label{eq:newF}
%\end{equation}
Then, the dynamics can be written as
\begin{equation}
  \dot x =F(x,\mu), \ \mu\in U(x)
  \label{eq:newdynamics}
\end{equation}
where 
\begin{equation}
  U(x)=\left\{\mu\Bigg|\ \|F(x,\mu)\|=\min\limits_{\stackrel{\nu_j\geq 0, j\in
  J(x)}{\nu_j=0,
  j\notin J(x)}}\|F(x,\nu)\|\right\}
  \label{eq:U}
\end{equation}
and $J(x)=\{j|\ h_j(x)\leq 0\}$. This formulation is given in the seminal work
by Rosen~\cite{Rosen:1965tg} along with the theorem that states that for any
initial condition in $\mc{C}$, a continuous solution $x(t)$ to \eqref{eq:newdynamics}
exists such that $x(t)\in \mc{C}$ for all $t>0$. Thus, we have the following
results. 
\begin{proposition}[Theorem 8~\cite{Rosen:1965tg}]
  The dynamical system \eqref{eq:newdynamics} is asymptotically stable
      on $\mc{C}$ if $D\omega(x, \mu)$ has eigenvalues in the open left-half plane
      for $x\in \mc{C}$ and $\mu\in U(x)$.
    \end{proposition}
    Further, if $x^\ast\in \mc{C}$ is a differential Nash equilibrium, we can linearize
    $\omega$ around $x^\ast$ and get the following sufficient condition
    guaranteeing $x^\ast$ attracts nearby strategies under the gradient flow
    $F(x,\mu)$.
    \begin{proposition}
      If $x^\ast\in \mc{C}$ is a differential Nash equilibrium, 
      %i.e. $\omega(x^\ast)=0$ and
  %$D_{ii}\ f_i(x^\ast)< 0$, 
  and the
  eigenvalues of $D\omega(x^\ast, \mu^\ast)$ are in the open left-half plane, then
  $x^\ast$ is an exponentially stable fixed point of the continuous-time
  dynamical system \eqref{eq:dynamics}.
\end{proposition}
Note that since in our estimation, we restrict $\theta_i\geq 0$, the $f_i$ will
be concave; hence, Nash equilibria of the game will be differential Nash
equilibria. 

These results imply that we can simulate the dynamical system in
\eqref{eq:newdynamics} in order to compute Nash equilibria of the game. Using a forward Euler
discretization scheme and a subgradient projection method, we can compute Nash
equilibria of the constrained game. The subgradient projection method is known
to converge to the unique Nash equilibrium of the constrained $n$-person concave game~\cite{flam:1990aa}.

%We remark that since the occupants are utility maximizers, the conditions for
%differential Nash equilibria as they appear in 
%The occupants vote by selecting their lighting
%preference. 
%%%%%%%%%%%%%%%%%%%%%%%%%%%%%%%%%%%%%%%%%%%%%%%%%%%%%%%%%%%%%%%%%%%%%%%%%%%%%%%
%%%%%%%%%%%%%%%%%%%%%%%%%%%%%%%%%%%%%%%%%%%%%%%%%%%%%%%%%%%%%%%%%%%%%%%%%%%%%%%
%%%%%%%%%%%%%%%%%%%%%%%%%%%%%%%%%%%%%%%%%%%%%%%%%%%%%%%%%%%%%%%%%%%%%%%%%%%%%%%

%% file: experiment.tex
% experimental setup file
The social game for energy savings that we have designed is such that 
occupants in an office builing vote according to their usage preferences of shared resources and are
rewarded with points based on how \emph{energy efficient} their strategy is in
comparison with the other occupants.  Having points increases the likelihood of
the occupant winning in a lottery. The prizes in the lottery consist of three Amazon
gift cards.

We have installed a
Lutron\footnote{\tt http://www.lutron.com/en-US/Pages/default.aspx} system for the control of the lights in the office. 
This system allows us to precisely control the lighting level of each of the
lights in the office. We use it to set the default lighting level as well as
implement the average of the votes each time the occupants change their lighting
preferences.

\begin{figure}[h]
  \begin{center}
    \subfigure[]{\includegraphics[scale=0.25]{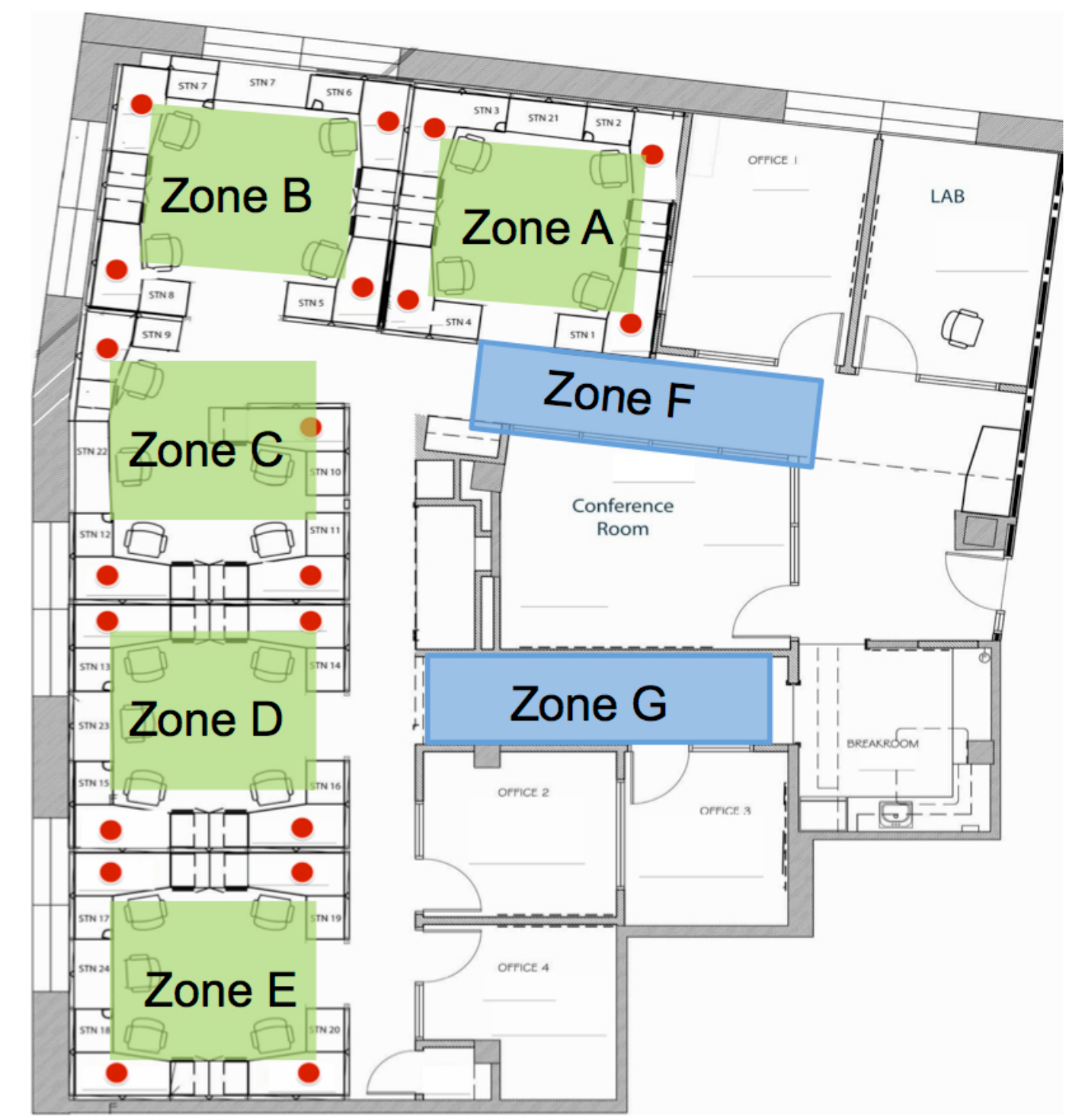}\label{fig:zones2}}
    \subfigure[]{\includegraphics[scale=0.1]{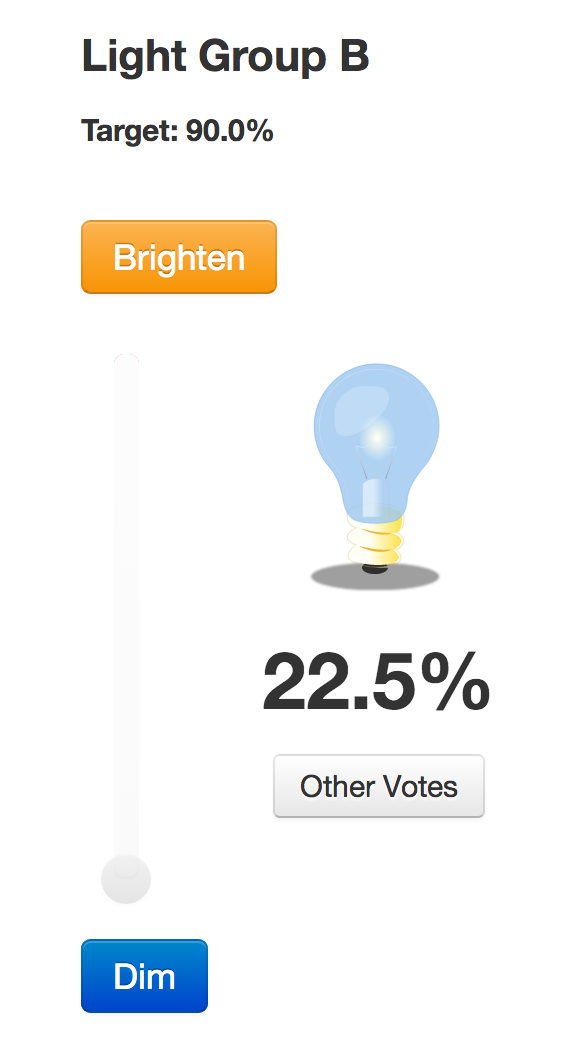}\label{fig:votes}}
    \subfigure[]{\includegraphics[scale=0.35]{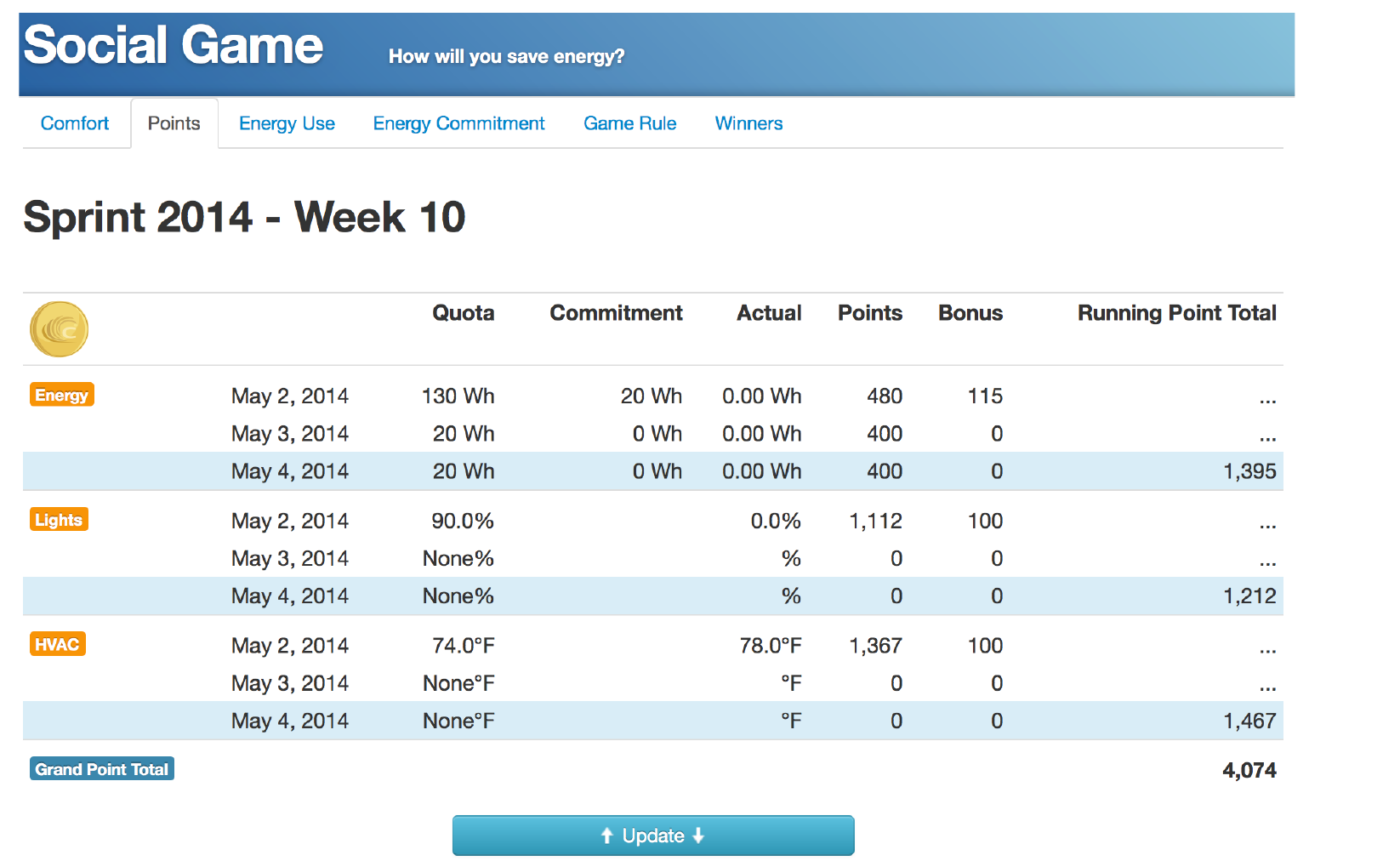}\label{fig:points}}
  \end{center}
  \caption{(a) Display of HVAC and lighting zones. Zones A-B are the five lighting
  zones and zones F-G are the two HVAC zones. (b) Display of how occupants can log their lighting
  vote. (c) Display of an occupant's point balance.}
  \label{fig:all}
\end{figure}

We have divided the office into five lighting zones and 
each zone has four occupants. Thus, there are 20 occupants who participate in
the social game. In addition, we have two heating, ventalating and air
conditioning (HVAC) zones and each zone has ten occupants (see
figure~\ref{fig:zones2}). 

We have developed an online platform in which the occupants can
login and participate in the game. This includes the ability for the occupants
to vote on their lighting and heating, ventilating and air conditioning (HVAC)
preferences as well as view all occupant point balances and all occupant
consumption patterns including the ability to monitor individual occupant
plug-load consumption. Figure~\ref{fig:votes}
shows a display of how an occupant can select their lighting preference and
Figure~\ref{fig:points} shows a sample of how occupants can see their point balance. 
%\begin{figure}[h]
%  \begin{center}
%%    \subfigure[]{\includegraphics[scale=0.1]{figs/visuallight.png}\label{fig:lightdis}}
%%\subfigure[]{\includegraphics[scale=0.2]{figs/zone.pdf}
%% \label{fig:zones}}
%\includegraphics[scale=0.1]{figs/visuallight.png}\includegraphics[scale=0.2]{figs/zone.pdf}
%      \end{center}
%      \label{fig:zones}
%      \caption{On the left is a visualization of the office floor plan and the lights above the cubicles where
%  dark green indicates a lower lighting setting and light green a higher
%  setting. On the right is a display of HVAC and lighting zones. Zones A-B are the five lighting
%  zones and zones F-G are the two HVAC zones.}
%\end{figure}

In this paper, we focus on a game focused on encouraging occupants to select
lower lighting settings in exchange for a chance to win in a lottery.  An occupant's vote is for the lighting level in their zone as well as for
neighboring zones. 
 The lighting setting that is implemented is the average of all the votes. 

Each day when an occupant logs into the online platform the first time after
they enter the office, they are considered present for the remainder of the day.
%Upon logging in the occupant can log their votes, view their point balance as
%well as other occupants' points, and view all occupant past consumption
%patterns. 
There is a default lighting setting. An occupant can leave the
lighting setting as the default after logging in or they can change it to some
other value in the interval $[0,100]$ depending on their preferences.

Some of the energy savings we achieve is due
to the default setting and some due to the social game. We are currently
conducting experiments to determine how much savings is due to the social game.
It is the building
managers duty to ensure that the occupants are satisfied (via appropriate lighting
level) and the building is operating in an energy efficient way. We believe that through optimal design of the incentives, we will be
able to achieve greater energy savings than would be possible by only adjusting
the default lighting setting. We leave this for future work.

%\begin{figure}[h]
%  \begin{center}
%    \includegraphics[scale=0.125]{figs/light.png}\includegraphics[scale=0.45]{figs/points.pdf}
%\label{fig:pointlight}
%  \end{center}
%  \caption{On the left is a display of how occupants can log their lighting
%  vote. On the right is a display of an occupant's point balance.}
%  %where the colors indicate how energy
%  %efficient occupants are currently on a scale from red (indicating high
%  %consumption) to green (indicating low consumption).}
% % \label{fig:light}
%\end{figure}

%The control of the lights' dim level as well as of the HVAC has been done with
%the implementation of the Lutron System (Lillian, can we give the website of the
%system that e use. It is a commercial one and I don't know what we have to include. the portal is: http://www.lutron.com/en-US/Pages/default.aspx ) via BACnet.

%This platform includes the capability for occupants to select both lighting and
%HVAC resources as well as monitor
%their personal plug-load consumption; however, the game we focus on in this
%paper is centered around the occupant lighting choices. 

%In particular, occupants login and select their lighting choice for the lights
%in the office any time they are present in the office (see
%figure~\ref{fig:pointlight}). 

%% file: results.tex
In this section, we report the results on the savings acheived through the game,
the utility learning problem as well as simulation of the estimated utilities. 

We use the data collected over the period from Mar. 3, 2014 to Jun. 5, 2014 when
occupants have regular working schedules in the office. The baseline lighting,
$x_b$, is 90\%, which is the standard lighting level prior to the beginning of
the experiment. Throughout this period, we have changed the default lighting
level three times (see Table~\ref{tab:default}).
\begin{table}
  \centering
  \begin{tabular}{|c|c|}
    \hline
    Period & Default Level \\
    \hline\hline
    March 3--April 10 & 20 \%\\
    April 11--May 1 & 10 \%\\
    May 2--May 23 & 60 \%\\
    May 24 -- June 5 & 90 \%\\
    \hline
  \end{tabular}
  \caption{Default levels for four periods during the experiment. By changing
  the default setting to $90$\% we isolate the savings due to the social game
  from those achieved by changing the default setting.}
  \label{tab:default}
\end{table}
%first, from 20 to 10 on Apr. 11, then from 10 to 60 on May 2, and lastly from 60 to 90 on May 24. 
%This offers the opportunity to separate the effects of social game on energy saving from that achieved by simply lower the default lighting. 
We divide each day into four regions based on the outside lighting in Berkeley during the summer, namely from 5 to 10am (Dawn), 10am to 5pm (Daylight), 5pm to 8pm (Dusk), and 8pm to the next day 5am (Night). The data is further processed by taking the average of votes in each region of the day for each user.

\subsection{Savings}
\begin{figure}[h]
  \begin{center}
    \includegraphics[scale=0.4]{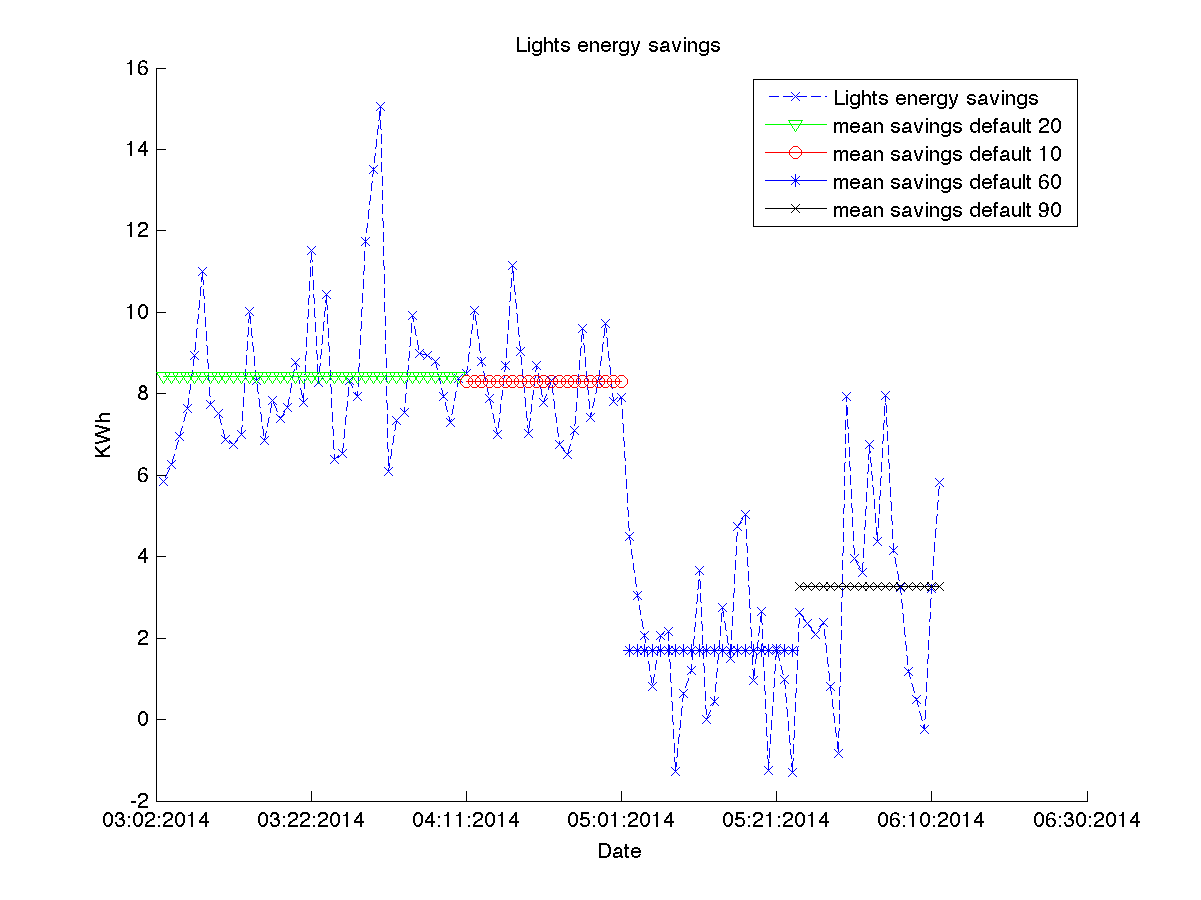}
  \end{center}%\vspace{-0.2in}
  \caption{Savings acheived per day (KWh). The mean savings over the four
  periods in which the default varied. Notice that in the period during which
  the default setting was at the baseline, there is still a savings of around
  3KWh per day.}
  \label{fig:savings}
\end{figure}

First, we highlight the savings achieved as a result of instituting the
social game. In Figure~\ref{fig:savings}, we report the savings per day in
KWh for the four periods in which the default varied. We remark that in the
last period in which the default setting was set to 90\% (which is the baseline
line setting), we still acheived a savings of 3KWh on average per day. Using the
mean savings for each of the periods and a rate of \$0.12/KWh, we estimate that
we saved \$73. In addition, over the period of 101 days
that the experiment was conducted the office consumed 2,185 KWh for lighting and
we saved approximately 601KWh. That is a 27.5\% reduction in energy. This savings is just due to a change in lighting usage behavior
for one small portion of a building. 

Our platform has the capability of including HVAC and plug-load in addition to
lighting. We plan to implement a similar social game in Singapore and we expect
much greater savings. This current experiment shows that a social game is a
viable way to engage building occupants and induce behavioral change toward more
energy efficient behaviors.

\subsection{Estimation}
%%%%%%%%%%%%%%%%%%%%%%%%%%%%%%%%%%%%%%%%%%%%%%%%%%%%%%%%%%%%%%%%%%%%%%%%%%%%%%%
The estimation proposed is Section 2.1 is performed for each user in each day
interval and default lighting interval. Only true votes, not the default votes,
are considered. We apply the bootstrapping method to obtain the empirical
distribution of $\theta_i$ for $i\in \{1,...,20\}$ by randomly sampling a subset
from the data \cite{efron:1994aa}. The mean and standard deviation for the users
which are the most active are reported in Table~\ref{tab:estimation}.

\begin{table}
\centering
\caption{Estimated utility parameter for selected set of active users. A, B, C,
D stand for the periods Dawn, Daylight, Dusk, and Night respectively. The
standard deviation is indicated inside the parentheses and the mean is given
outside of the parentheses. Blank indicates that the occupant did not vote during that
period. Hence, they have no estimated utility. Occupants' whose number is in
boldface have won at least once in the lottery. Note that numbers marked with * indicate that we do not have reliable estimates due to limited data.}
    \begin{tabular}{cc|c|c|c|c|c|c|}
\cline{3-8}
& & \multicolumn{6}{ |c| }{Active users (selected)} \\ \cline{3-8}
& & 2 & {\bf 6} & {\bf 8} & 10 &{\bf 14} & 20  \\ 
\cline{1-8}
\multicolumn{1}{ |c| }{\multirow{4}{*}{Default 20}}  &
\multicolumn{1}{ |c| }{A} & .26(.25) & 4657(481) & 3671(126) & 4658(0*) & 3054(141)& 3857(110)     \\ \cline{2-8}
\multicolumn{1}{|c|}{} &
\multicolumn{1}{ |c| }{B} & .41(.09) & 2932(99) & 3386(73) & 271(307)& 3350(96) & 691(473)      \\ \cline{2-8}
\multicolumn{1}{|c|}{} &
\multicolumn{1}{ |c| }{C} & .00(.01) & 1808(1008) & 3290(194) & 1220(0*)   & 3332(164) & 1222(259)   \\ \cline{2-8}
\multicolumn{1}{|c|}{} &
\multicolumn{1}{ |c| }{D} & .00(.00) & 822(1465) & 3700(461) & 3420(295)  & 3756(575) & 1095(285)     \\ \cline{1-8}

\multicolumn{1}{ |c| }{\multirow{4}{*}{Default 10}}  &
\multicolumn{1}{ |c| }{A} & .96(.39) & 294(759) & 2923(215) & 2446(0*) & 2971(508)& 195(335)     \\ \cline{2-8}
\multicolumn{1}{|c|}{} &
\multicolumn{1}{ |c| }{B} & .24(.60) & 833(796) & 2847(320) & 2042(0*)& 3219(339) & 258(339)      \\ \cline{2-8}
\multicolumn{1}{|c|}{} &
\multicolumn{1}{ |c| }{C} & .07(.12) & 0(0*) & 2924(224) & 3485(0*)   & 670(441) & 643(469)   \\ \cline{2-8}
\multicolumn{1}{|c|}{} &
\multicolumn{1}{ |c| }{D} & .09(.25) & 625(816) & 3542(474) & 3305(0*)  & 1793(1187) & 824(534)     \\ \cline{1-8}

\multicolumn{1}{ |c| }{\multirow{4}{*}{Default 60}}  &
\multicolumn{1}{ |c| }{A} & .28(.59) & 469(1717) & 6790(1267) &   & 3180(827)& 504(940)     \\ \cline{2-8}
\multicolumn{1}{|c|}{} &
\multicolumn{1}{ |c| }{B} & .07(.19) & 1062(1135) & 5741(734) & 6327(199)& 6180(881) & 104(484)      \\ \cline{2-8}
\multicolumn{1}{|c|}{} &
\multicolumn{1}{ |c| }{C} & .00(.00) & 1146(1927) & 6166(502) & 3752(0*)   & 7856(1728) & 588(903)   \\ \cline{2-8}
\multicolumn{1}{|c|}{} &
\multicolumn{1}{ |c| }{D} & .12(.18) & 3947(2434) & 6670(0*) & 5296(0*)  & 3628(3394) & 881(4)     \\ \cline{1-8}

\multicolumn{1}{ |c| }{\multirow{4}{*}{Default 90}}  &
\multicolumn{1}{ |c| }{A} & .01(.01) & 9045(1562) & 7835(2465) &   &  & 3333(0*)     \\ \cline{2-8}
\multicolumn{1}{|c|}{} &
\multicolumn{1}{ |c| }{B} & .00(.01) & 7624(1699) & 9479(926) &  &   & 1923(2010)      \\ \cline{2-8}
\multicolumn{1}{|c|}{} &
\multicolumn{1}{ |c| }{C} & .02(.03) & 8962(947) & 8761(983) &     &   & 3333(0*)   \\ \cline{2-8}
\multicolumn{1}{|c|}{} &
\multicolumn{1}{ |c| }{D} &   &   & 5000(461) & 5000(0*)  &   &       \\ \cline{1-8}
\end{tabular}
  \label{tab:estimation}
\end{table}

\begin{figure}[h]
  \begin{center}
    {\includegraphics[scale=0.45,trim = 0 15mm 0 0]{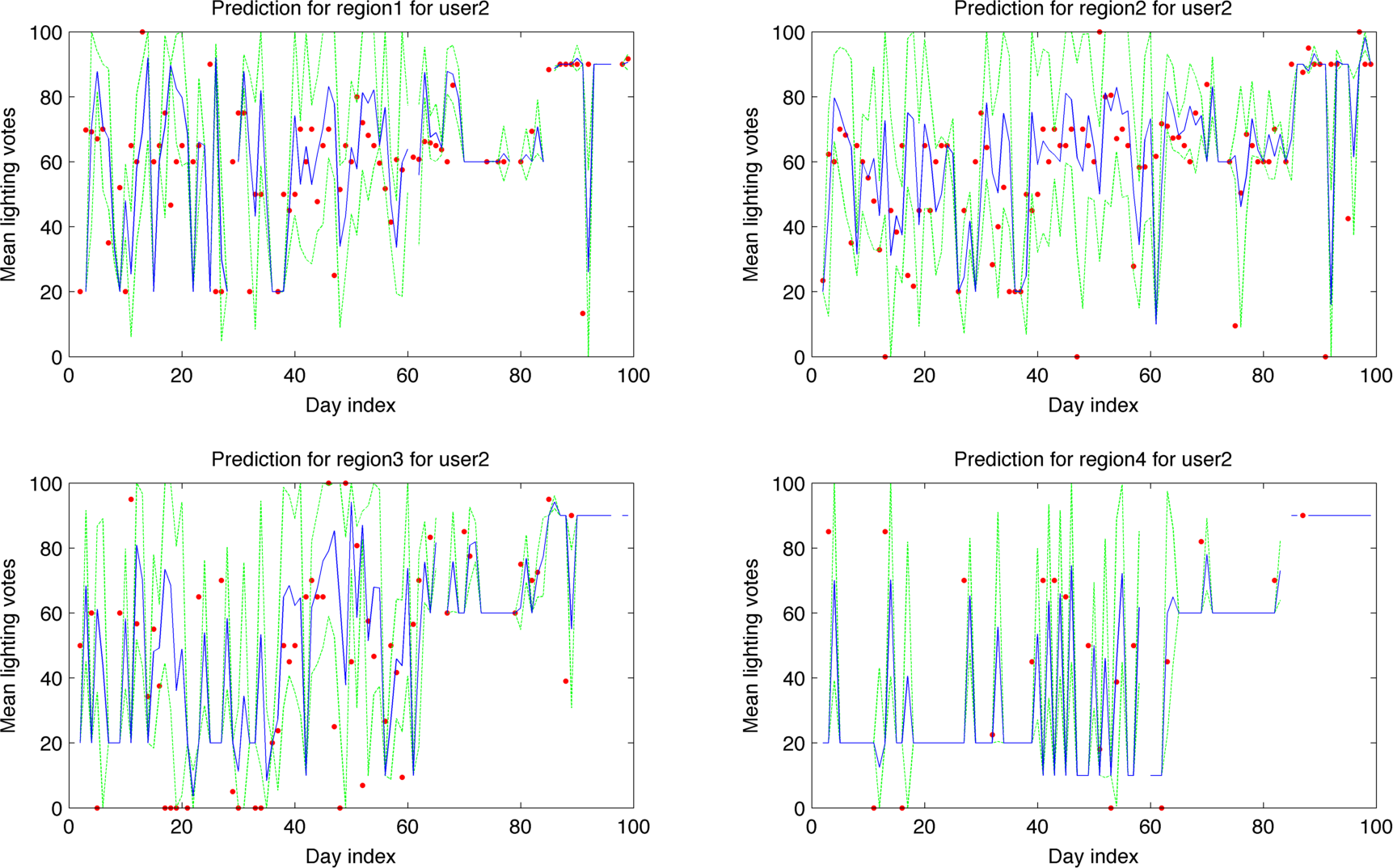}}
  \end{center}
  \vspace{0.1in}
  \caption{One day ahead prediction by the Nash equilibrium algorithm for
  occupant 2. Red: ground truth, i.e. occupant 2's actual votes. We sample
  from the distributions across the events \emph{absent}, \emph{active},
  \emph{default} for each occupant and simulate the game with the set of
  \emph{active} and \emph{default} players. We repeat
  this 20 times for each day and generate a distribution for the predictions of each
  occupant's behavior. Blue: mean of
  prediction. Green: prediction within 1 standard deviation of the prediction
  mean. Gaps in the plots indicate that the occupant was not present on that
  day.}
\label{fig:simulationuser2}
\end{figure}
\begin{figure}[h]
  \begin{center}
    {\includegraphics[scale=0.45,trim = 0 15mm 0 0]{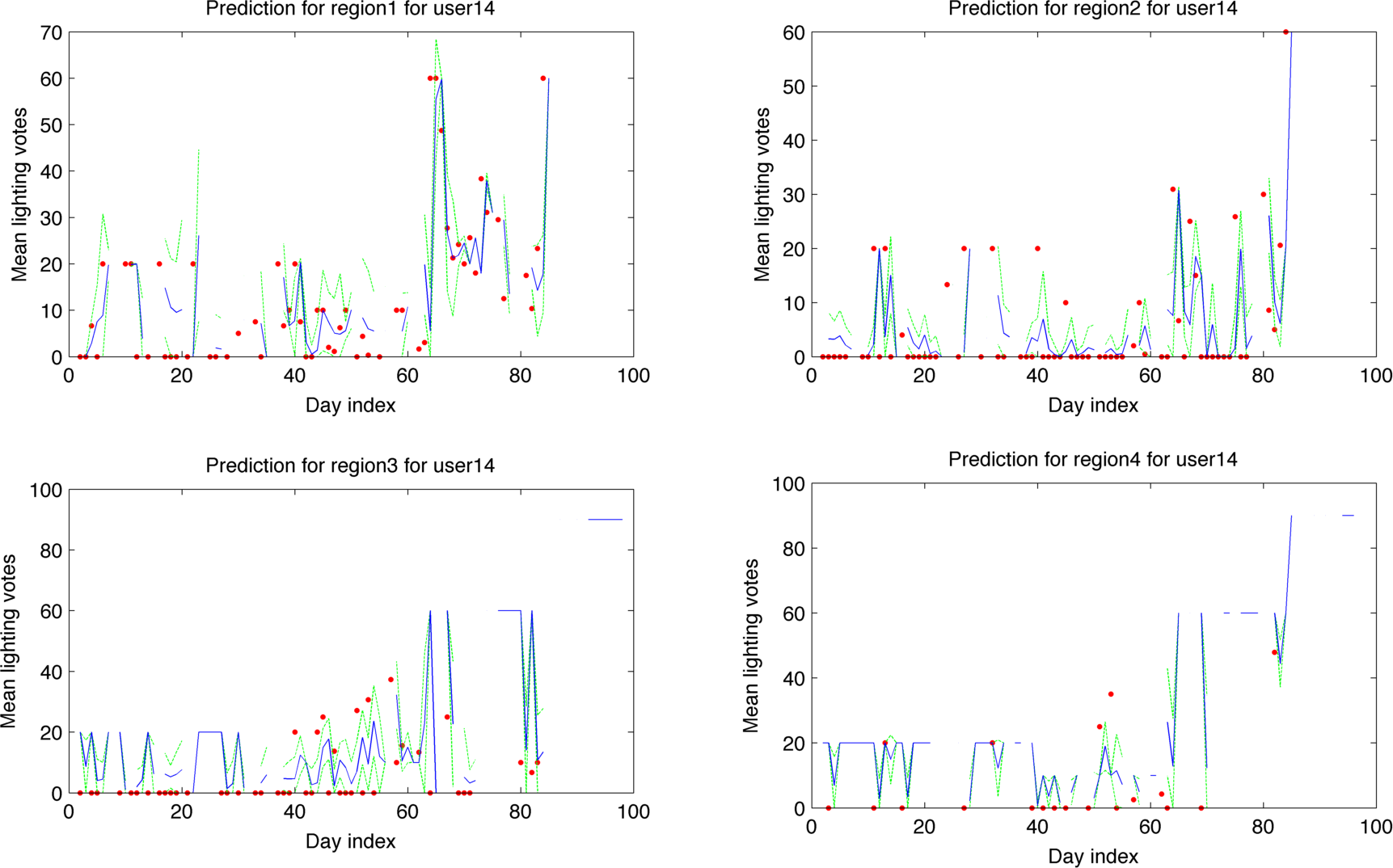}}
  \end{center}
  \vspace{0.1in}
  \caption{One day ahead prediction by the Nash equilibrium algorithm for
  occupant 14. Red: ground truth, i.e. occupant 14's actual votes. We sample
  from the distributions across the events \emph{absent}, \emph{active},
  \emph{default} for each occupant and simulate the game with the set of
  \emph{active} and \emph{default} players. We repeat
  this 20 times for each day and generate a distribution for the predictions of each
  occupant's behavior. Blue: mean of
  prediction. Green: prediction within 1 standard deviation of the prediction
  mean. Gaps in the plots indicate that the occupant was not present on that
  day.}
\label{fig:simulationuser14}
\end{figure}
We remark that occupant 2 has a very low mean for the parameter $\theta_2$ as
compared to the other active occupants.
By examining the ground truth values (red dots) in
Figure~\ref{fig:simulationuser2}, we see that occupant 2 often votes for a
lighting setting around 60-70\%. On the other hand, in
Figure~\ref{fig:simulationuser14}, we can see the ground truth of occupant 14 who
often votes for a lighting setting of 0\%. This player is more aggressive than
occupant 2 and this behavior is reflected in the mean of the estimate
for parameter $\theta_{14}$.

\subsection{Simulation}
\begin{figure}[h]
  \begin{center}
    {\includegraphics[scale=0.45,trim = 0 15mm 0 0]{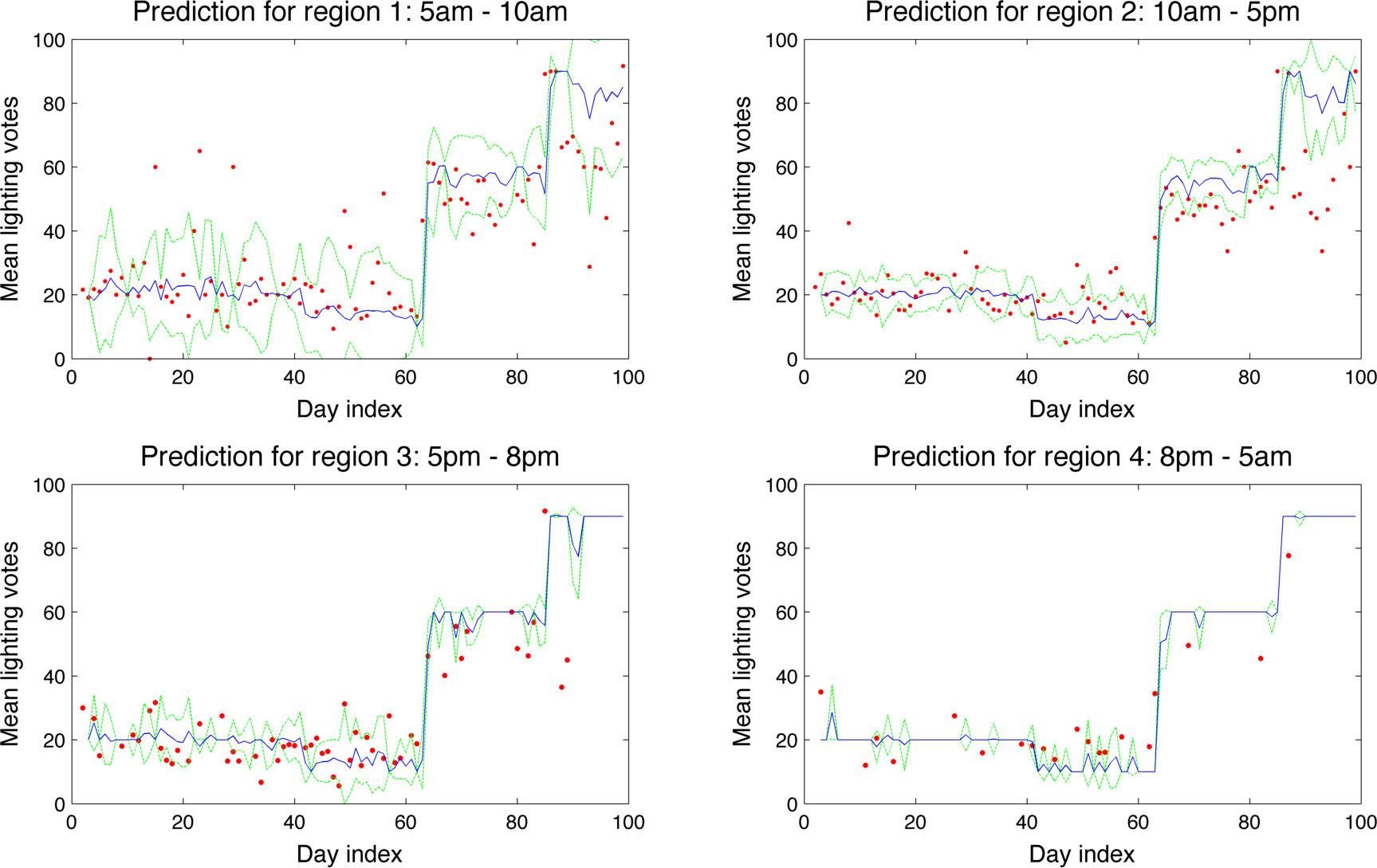}}
  \end{center}
  \vspace{0.1in}
  \caption{One day ahead prediction by the Nash equilibrium algorithm for the
   average implemented lighting setting in the office each day. We sampled from
   the distribution over the events \emph{absent}, \emph{active}, \emph{default}
   for each of the players and simulated Nash given the sample of \emph{active}
   and \emph{default} players. We
   repeated this 20 times for each day and generated a distribution for the prediction of the
   implemented lighting setting. Red:
  ground truth, i.e. the average lighting setting that is implemented per day.
  Blue: mean of prediction as given by the Nash simulation. Green: prediction within 1 standard
  deviation of the prediction mean.}
\label{fig:simulation}
\end{figure}

%%%%%%%%%%%%%%%%%%%%%%%%%%%%%%%%%%%%%%%%%%%%%%%%%%%%%%%%%%%%%%%%%%%%%%%%%%%%%%%
%%%%%%%%%%%%%%%%%%%%%%%%%%%%%%%%%%%%%%%%%%%%%%%%%%%%%%%%%%%%%%%%%%%%%%%%%%%%%%%
%%%%%%%%%%%%%%%%%%%%%%%%%%%%%%%%%%%%%%%%%%%%%%%%%%%%%%%%%%%%%%%%%%%%%%%%%%%%%%%

To capture the working schedules of each user, we employ a simple probabilistic model which determines the probability of individual user being absent, $p_i^{\text{absent}}$, present and playing default, $p_i^{\text{present, default}}$, and present and actively playing, $p_i^{\text{present, active}}$. By assumption the sample space $\Omega$ includes the above three outcomes, and the probability mass functions should sum to unity. This probability is estimated by $p_i^E=\tfrac{N_{i,E}}{N_i}$, where $E$ is the event of one of the three outcomes, $N_{i,E}$ is the number of event $E$ for user $i$, and $N_i$ is the number of total events.

For the prediction of the next day lighting votes, we randomly sample from this
distribution to determine the set of active, default, and absent users, then
obtain a local Nash equilibrium for them. This step is performed 20 times for each day to predict the distribution of votes, as shown in Figure~\ref{fig:simulation}.

As can be seen, the Nash equilibrium captures substantial variations in
the data. We also compared the results of prediction with the autoregressive
integrated moving average (ARIMA) model~\cite{hannan:2012aa}, constant model which uses the default
lighting for prediction, and the persistent model which uses the previous day
value for prediction. The mean squared errors (MSE) of the models are summarized
in Table~\ref{tab:prediction}. The Nash equilibrium achieves a prediction that
is the most accurate as compared with other models, which presents it favorably
for leaders in the Stackelberg game to design optimal incentives to motivate
energy saving behaviors. Indeed in the Stackelberg framework, the leader (building
manager) assumes that the agents (occupants) are utility maximizers and play
Nash. Hence, we will be able to integrate our estimation algorithm into an
online algorithm for designing incentives.  

\begin{table}
\centering
\caption{Mean square error (MSE) of four algorithms that predict the one day
ahead occupant behavior over the period of study (101 days): ARIMA(1,0,1) (we use one
autoregressive term, zero nonseasonal differences, and one lagged forecast error
in the prediction equation)~\cite{hannan:2012aa}, Nash, a model which uses
the default as the prediction, and a model which uses occupants' previous votes
as the prediction.  Nash out performs each of the other methods.}
    \begin{tabular}{| l | c | c | c |c|}
    \hline
   \textit{ Model } & ARIMA(1,0,1) & Nash & Constant & Persistent \\ \hline
     \textit{ MSE }& 13.92 & \textbf{12.46 }& 16.96 & 13.42\\
    \hline
    \end{tabular}
    
  \label{tab:prediction}
\end{table}

%% file: discussion.tex
%\begin{itemize}
%  \item should also discuss environmental noise model and where we are going
%    with that.
%\end{itemize}
We have designed and implemented a social game for inducing building occupants
to behave in an energy efficient manner. We presented data and results
pertaining to the game in which occupants select their lighting preferences and
win points depending on how far their vote is from the baseline lighting setting
and proportional to other occupants' votes distances from the baseline. As a
result, the occupants are interacting in a competitive environment which we model as a
non-cooperative game. We show that we get significant savings as compared to
usage prior to the implementation of the social game. This savings is due to
both a change in the default setting as well as due to the incentives offered in
the social game.

We described the experimental set-up which includes an online platform for the
implementation of the social game as well as the use of a Lutron lighting
system for percise control of the lighting setting.
Our platform also includes the ability to implement a
social game centered around HVAC settings as well as occupant plug-load
consumption. We leave exploring these additional features as future work. 

We have formulated the problem of estimating the occupant utility functions as a convex optimization
problem and estimated occupant utilities in a 20 player social game. We simulated the game using the
estimated utility functions and showed that our model is a good predictor for
occupant behavior. It out performs a number of other estimation techniques
including ARIMA. 

There are several ways in which we believe we can improve our estimate of
the utility functions of the occupants. We did not consider the environmental
noise such as variations in natural light. We instead used a heuristic to capture this
variation by breaking the day into intervals in which the natural light entering
the office is most consistent. In addition, we did not consider any
information on the occupants' schedules or location in the office with respect to
windows. We could incorporte these aspects into our estimation as priors
on the parameters of the occupants utility function or as
a noise process in the estimated behavior model. We
leave this as future work.

In the experiments used for this paper, we selected the value of
$\rho$ based on heuristics. Our goal is to design $\rho$ in an optimal way.
We can leverage the fact that we have modeled occupants as utility maximizers
who play in a non-cooperative game by considering the design of $\rho$ by the
building manager. In particular, we can model this interaction between the
building manager and the occupants as a Stackelberg game. In this framework, the
building manager would perform an online estimation of the occupants' utility
functions and update $\rho$ accordingly~\cite{ratliff:2014aa}. We believe that
by optimizing the incentive $\rho$, we can achieve greater savings. We are
currently implementing such a scheme in our experimental platform.